\documentclass[11pt]{amsart}

\usepackage{latexsym}
\newtheorem{thm}{Theorem}[section]
\newtheorem{lem}{Lemma}[section]
\newtheorem{prop}{Proposition}[section]

\newtheorem{rmk}{Remark}[section]

\def\1{{{\mbox{${\rm{1\negthinspace\negthinspace I}}$}}}}

\newcommand{\E}{\mathbb{E}}
\newcommand{\R}{\mathbb{R}}

\begin{document}
\title[Nonparametric estimation for L\'evy processes]{Nonparametric adaptive estimation for pure jump L\'evy processes.}
\author[F. Comte]{F. Comte}
\author[V. Genon-Catalot]{V. Genon-Catalot$^{1}$}\thanks{$^1$ Universit\'e Paris V,  MAP5, UMR CNRS
8145. email: fabienne.comte@univ-paris5.fr and genon@math-info.univ-paris5.fr}

\begin{abstract} This paper is concerned with nonparametric estimation of the L\'evy density of a pure jump L\'evy process. The sample path is observed at $n$ discrete instants with fixed sampling interval. We construct a collection of estimators obtained by deconvolution methods and deduced from appropriate estimators of the characteristic function and its first derivative. We obtain a bound for the ${\mathbb L}^2$-risk, under general assumptions on the model. Then we propose a penalty function that allows to build an adaptive estimator. The risk bound for the adaptive estimator is obtained under additional assumptions on the L\'evy density. Examples of models fitting in our framework are described and rates of convergence of the estimator are discussed. 
{\bf \today}
\end{abstract}
\maketitle

\noindent {\bf {\sc Keywords.}} Adaptive Estimation; Deconvolution; L\'evy process;
Nonparametric Projection Estimator.

\medskip\

\section{Introduction}
In recent years, the use of L\'evy processes for modelling purposes has become very popular in many areas and especially in the field of finance (see {\em e.g.} Eberlein and Keller (1995), Barndorff-Nielsen and Shephard (2001), Cont and Tankov (2004); see also Bertoin (1996) or Sato (1999) for a comprehensive study for these processes). The distribution of a L\'evy process is usually specified by its characteristic triple (drift, Gaussian component and L\'evy measure) rather than by the distribution of its independent increments. Indeed, the exact distribution of these increments is most often intractable or even has no closed form formula. For this reason, the standard parametric approach by likelihood methods is a difficult task and many authors have rather considered nonparametric methods. For L\'evy processes, estimating the L\'evy measure is of crucial importance since this measure specifies the jumps behavior. Nonparametric estimation of the L\'evy measure has been the subject of several recent contributions. The statistical approaches depend on the way observations are performed. For instance, Basawa and Brockwell (1982) consider non decreasing L\'evy processes and observations of jumps with size larger than some positive $\varepsilon$, or discrete observations with fixed sampling interval. They build nonparametric estimators of a distribution function linked with the L\'evy measure. More recently, Figueroa-L\'opez and Houdr\'e (2006) consider a continuous-time observation of a general L\'evy process and study penalized projection estimators of the L\'evy density based on integrals of functions with respect to the random Poisson measure associated with the jumps of the process. However,  their approach remains theoretical since these Poisson integrals are hardly accessible.  

In this paper, we consider nonparametric estimation of the L\'evy measure for real-valued L\'evy processes of pure jump type, {\em i.e.} without drift and Gaussian component. We rely on the common assumption that the L\'evy measure admits a density $n(x)$ on ${\mathbb R}$ and assume that the process is discretely observed with fixed sampling interval $\Delta$.  Let $(L_t)$ denote the underlying L\'evy process and  $(Z_k^{\Delta}= L_{k\Delta}- L_{(k-1)\Delta}, k=1, \ldots, n)$ be the observed random variables which are independent and identically distributed. Under our assumption, the characteristic function of $L_{\Delta}= Z_1^{\Delta}$ is given by the following simple formula:
\begin{equation} \label{fc}
\psi_{\Delta}(u)= {\mathbb E}(\exp{i u Z_1^{\Delta}})=\exp{(\Delta \int_{{\mathbb R}}(e^{iux} -1) n(x) dx)}
\end{equation}
where  the unknown function  is the  L\'evy density  $n(x)$.  It  is therefore
natural to investigate the  nonparametric estimation of $n(x)$ using empirical
estimators  of  the charasteristic  functions  and  its  derivatives and  then
recover the L\'evy density by  Fourier inversion. This approach is illustrated
by Watteel and  Kulperger (2003) and Neumann and  Reiss (2000). However, these
authors   consider  general   L\'evy  processes,   with  drift   and  Gaussian
component. Hence, at least two  derivatives of the characteristic function are
necessary to reach the L\'evy density.  Moreover, the way Fourier inversion is
done in concrete is not detailed in these papers. In our case, under the assumption that $\int_{{\mathbb R}} |x| n(x) dx < \infty$, we get the simple relation:
\begin{equation}\label{basic}
g^*(u)= \int e^{iux}g(x)dx=-i\frac{\psi'_{\Delta}(u)}{\Delta\psi_{\Delta}(u)},
\end{equation}
with $g(x)=x n(x)$. This equation indicates that we can estimate $g^*(u)$ by using empirical counterparts of $\psi_{\Delta}(u)$ and $\psi'_{\Delta}(u)$ only. Then, the problem of recovering an estimator of $g$ looks like a classical deconvolution problem. We have at hand the methods used for estimating unknown densities of random variables observed with additive independent noise. This requires the additional assumption that $g$ belongs to ${\mathbb L}^{2}({\mathbb R})$.  However, the problem of deconvolution set by equation (\ref{basic}) is not standard and looks more like deconvolution in presence of unknown errors densities. This is due  to the fact that both the numerator and the denominator are unknown and have to be estimated from the same data. This is why our estimator of 
$\psi_{\Delta}(u)$ is not a simple empirical counterpart. Instead, we use a truncated version analogous to the one used in Neumann (1997) and Neumann and Reiss (2000).  

Below, we show how to adapt the deconvolution method described in Comte {\it et al.}~(2006). We consider an adequate sequence $(S_m, m=1, \ldots, m_n)$ of subspaces of ${\mathbb L}^{2}({\mathbb R})$ and build a collection of projection estimators $({\hat g}_m)$. Then using a penalization device, we select through a data-driven procedure the best estimator in the collection. We study the ${\mathbb L}^{2}$-risk of the resulting estimator under the asymptotic framework that $n$ tends to infinity. Although the sampling interval $\Delta$ is fixed, we keep it as much as possible in all formulae since the distributions of the observed random variables highly depend on $\Delta$. 

In Section 2, we give assumptions and some preliminary properties. Section 3 contains examples of models included in our framework. Section 4 describes the statistical strategy. We present the projection spaces and define the collection of estimators. Proposition 4.1 gives the upper bound for the risk of a projection estimator on a fixed projection space. This proposition guides the choice of the penalty function and allows to discuss the rates of convergence of the projection estimators. Afterwards,  we  introduce a theoretical penalty (depending on the unknown characteristic function $\psi_{\Delta}$) and study the risk bound of a false estimator (actually not an estimator) (Theorem 4.1). Then, we replace the theoretical penalty by an estimated counterpart and give the upper bound of the risk of the resulting penalized estimator (Theorem 4.2). Section 6 gives some conclusions and open problems. Proofs are gathered in Section 6. In the Appendix, a fondamental result used in our proofs is recalled.

\section{Framework and assumptions.}
Recall that we consider the discrete time observation with sample step $\Delta$ of a L\'evy process $L_t$ with L\'evy density $n$ and characteristic function given by (\ref{fc}). 
We assume that $(L_t)$ is a pure jump process with finite variation on compacts.
When the L\'evy measure $n(x)dx$ is concentrated on $(0, +\infty)$, then
$(L_t)$ has increasing paths and is called a subordinator.
We focus on the estimation of the real valued function
\begin{equation}
g(x)= x n(x),
\end{equation}
and introduce the following assumptions on the function $g$:
\begin{itemize}
\item[(H1)$\;\;\;\;\;$]  $\int_{\R} |x| n(x) dx < \infty.$
\item [(H2$(p)$)] For $p$ integer,  $\int_{\R}
  |x|^{p-1} |g(x)| dx <\infty$.
\item [(H3)$\;\;\;\;\;$] The function $g$ belongs to ${\mathbb L}_{2}(\R)$.
\end{itemize}
Note that (H1) is stronger than the usual assumption $\int (|x|\wedge 1) n(x)dx <+\infty$, and is also a moment assumption for $L_t$.
Under the usual assumption, (H2$(p)$) for $p\geq 1$ implies (H1) and (H2$(k)$) for $k\leq p$. 

Our estimation procedure is based
on the random variables \begin{equation}\label{kiezi} Z_i^{\Delta}= L_{i\Delta}-
L_{(i-1)\Delta} , i= 1, \ldots, n,\end{equation}  which are independent,
identically distributed, with common characteristic function
$\psi_{\Delta}(u)$.

The moments of $Z_1^{\Delta}$ are linked with the function $g$. More precisely,
we have:
\begin{prop} Let $p \ge 1$ integer. Under (H2)($p$),  $\E|Z_1^{\Delta}|^{p} <\infty$. Moreover, setting, for $k=1, \dots p$, $M_k= \int_{\R} x^{k-1}g(x)dx$, we have ${\mathbb E}(Z_1^{\Delta})=\Delta M_1$, ${\mathbb E}[(Z_1^{\Delta})^2]=\Delta M_2+\Delta^2 M_1$, and more generally, $\E[(Z_{1}^{\Delta})^{l}]= \Delta\; M_l + o(\Delta)$ for all $l=1,\ldots, p$.
\end{prop}

\noindent {\bf Proof.} By the assumption,  the exponent of the exponential in (\ref{fc}) is $p$
times differentiable and, by derivating $\psi_{\Delta}$, we get the result.
$\Box$\\

Assumption (H1) yields the relation (\ref{basic}), which is the basis of our
estimation procedure. 
We need a precise control of $\psi_{\Delta}$. For this, we introduce the assumption that, for $m_n$ an integer to be defined later, the following holds:
\begin{itemize} \item[(H4)$\;\;\;\;\;$]  $\forall x\in {\mathbb R}, \mbox{ we have }  c_\psi(1+x^2)^{-\Delta\beta/2} \leq |\psi_{\Delta}(x)|\leq C_\psi(1+x^2)^{-\Delta \beta/2},$
\end{itemize}  for some given constants $c_{\psi}, C_{\psi}$ and $\beta\geq 0$. Note that an assumption of this type is also considered in Neumann and Reiss~(2007).

For the adaptive version of our estimator, we need additional assumptions for $g$:
\begin{itemize} \item[(H5)$\;\;\;\;\;$]  There exists some positive $a$  such that $\int |g^*(x)|^2(1+x^2)^adx < +\infty$,
\end{itemize} and 
\begin{itemize} \item[(H6)$\;\;\;\;\;$] $\int x^2g^2(x)dx<+\infty.$
\end{itemize}
We must set independent assumptions for $\psi_{\Delta}$ and $g$, since there may be no relation at all between these two functions (see the examples). Note that, in Assumption (H5), which is a classical regularity assumption, the knowledge of $a$ is not required. 
 
\section{Examples.}\label{example}
\subsection{Compound Poisson processes.}
Let $L_t= \sum_{i=1}^{N_t} Y_i$, where $(N_t)$ is a Poisson process
with constant intensity $c$ and $(Y_i)$ is a sequence of i.i.d. random
variables with density $f$ independent of the process $(N_t)$. Then,
$(L_t)$ is a compound Poisson process with characteristic function
\begin{equation}
\psi_{t}(u)= \exp{c t \int_{\R}(e^{iux} -1) f(x) dx}.
\end{equation}
Its L\'evy density is $n(x)= c f(x)$. Assumptions (H1)-(H2)$(p)$ are
equivalent to $\E(|Y_1|^{p})<\infty$. Assumption (H3) is equivalent to
$\int_{\R} x^2 f^{2}(x) dx < \infty$, which holds for instance if $\sup_x f(x)<+\infty$ and ${\mathbb E}(Y_1^2)<+\infty$. We can compute the distribution
of $Z_1^{\Delta}=L_{\Delta}$ as follows:
\begin{equation}
  P_{Z_1^{\Delta}}(dz) = e^{-c \Delta} ( \delta_{0}(dz) + \sum_{n \ge 1} f^{*n}(z) \frac{(c \Delta)^{n}}{n!} dz ).
\end{equation}
We have the following bound:
\begin{equation}
1\geq |\psi_{\Delta}(u)| \ge e^{-2 c \Delta}.
\end{equation}
On this example, it appears clearly that we can not link the regularity assumption on $g$ and (H4) which holds with $\beta=0$.

\subsection{The L\'evy gamma process.}
Let $\alpha>0,\beta>0$. The L\'evy gamma process $(L_t)$ with parameters $(\beta,\alpha)$ is a subordinator such that, for all $t>0$, $L_t$ has distribution Gamma with parameters $(\beta t,\alpha)$, {\em i.e.} has density:
\begin{equation}
\frac{\alpha^{\beta t}}{\Gamma(\beta t)} x^{\beta t-1} e^{-\alpha x} 1_{x \ge 0}.
\end{equation}
The characteristic function of $Z_1^{\Delta}$ is equal to:
\begin{equation}
\psi_{\Delta}(u)= \left(\frac{\alpha}{\alpha-i u}\right)^{\beta \Delta}.
\end{equation}
The L\'evy density is $n(x)= \beta x^{-1} e^{-\alpha x}\1_{\{x>0\}}$ so that
$g(x)= \beta e^{-\alpha x}\1_{\{x>0\}}$ satisfies our assumptions. We have:
\begin{equation}
\frac{\psi_{\Delta}^{'}(u)}{\psi_{\Delta}(u)}= i \Delta \frac{\beta}{\alpha-i u}, \quad \quad
|\psi_{\Delta}(u)|= \frac{\alpha^{\beta \Delta}}{(\alpha^2+u ^2)^{\beta\Delta/2}}
\end{equation}

\subsection{Another class of subordinators}\label{other}
Consider the L\'evy process $(L_t)$ with L\'evy density
$$
n(x)= c x^{\delta-1/2} x^{-1} e^{-\beta x} 1_{x>0}, 
$$
where $(\delta,\beta,c)$ are positive parameters.
If $\delta> 1/2$, $\int_0^{+\infty}n(x)dx<+\infty$, and we recover compound Poisson processes.
If $0<\delta\leq 1/2$, $\int_0^{+\infty}n(x)dx=+\infty$ and $g(x)=xn(x)$ belongs to ${\mathbb L}^2({\mathbb R})\cap
{\mathbb L}^1({\mathbb R})$.
The case $\delta=0$, which corresponds to the L\'evy inverse Gaussian process does not fit in our framework.
For $0<\delta< 1/2$, we find 
$$g^*(x)=c\frac{\Gamma(\delta+1/2)}{(\beta-ix)^{\delta+1/2}},$$ and $$|\psi_{\Delta}(x)|=\exp\left( -c\frac{\Delta\Gamma(\delta+1/2)}{1/2-\delta} [(\beta^2+x^2)^{-(\delta-1/2)/2}-\beta^{-(\delta-1/2)}]\right).$$
It is important to mention that $\psi_{\Delta}$ above does not satisfy assumption (H4) since \begin{equation}\label{equiv} |\psi_{\Delta}(x)|\sim_{x\rightarrow +\infty} K(\beta,\delta)\exp(-c\Delta\frac{\Gamma(\delta+1/2)}
{1/2-\delta}x^{-\delta+1/2})\end{equation}
where $K(\beta,\delta)=\exp\left( c\frac{\Delta\Gamma(\delta+1/2)}{1/2-\delta}\beta^{-(\delta-1/2)}\right)$.
Thus, it has an exponential rate of decrease.

\subsection{The bilateral Gamma process.}
This process has been recently introduced by K\"{u}chler and
Tappe~(2008).  Consider $X,Y$ two independent random variables,
$X$ with distribution $\Gamma(\beta,\alpha )$ and $Y$ with distribution
$\Gamma(\beta',\alpha')$. Then, $Z=X-Y$ has distribution bilateral gamma
with parameters $(\beta,\alpha,\beta',\alpha')$, that we denote by
$\Gamma(\beta,\alpha;\beta',\alpha')$. The characteristic function of $Z$ is equal
to:
\begin{equation}
\psi(u)= \left(\frac{\alpha}{\alpha-i u}\right)^{\beta} \left(\frac{\alpha'}{\alpha'+i
u}\right)^{\beta'}= \exp{( \int_{\R}(e^{iux} -1) n(x) dx)},
\end{equation}
with
$$
n(x)= x^{-1} g(x),
$$
and, for $x \in \R$,
$$
g(x)= \beta e^{-\alpha x} 1_{(0, +\infty)(x)} - \beta' e^{-\alpha'|x|} 1_{(-\infty,0)}(x).
$$
The bilateral Gamma process $(L_t)$ has characteristic function
$\psi_t(u)= \psi(u)^{t}$.

The method can be generalized and we may consider L\'evy processes on ${\mathbb R}$ obtained by bilateralisation of two subordinators.

\subsection{Subordinated Processes.}
 Let $(W_t)$  be a Brownian  motion, and let  $(Z_t)$ be an  increasing L\'evy
 process  (subordinator), independent  of  $(W_t)$. Assume  that the  observed
 process is $$L_t=W_{Z_t}.$$ We have $$\psi_{\Delta}(u)={\mathbb E}(e^{iuL_{\Delta}}) = {\mathbb E}(e^{-\frac{u^2}2 Z_{\Delta}}).$$ As $Z_t$ is 
positive, we consider, for $\lambda \ge 0$,
$$\vartheta_{\Delta}(\lambda)   =  {\mathbb   E}(e^{-\lambda   Z_{\Delta}})  =
\exp\left( -\Delta  \int_0^{+\infty}(1-e^{-\lambda x})n_Z(x)dx\right),$$ where
$n_{Z}$  denotes  the  L\'evy density  of  $(Z_t)$.  Now  let us  assume  that
$g_Z(x)=xn_Z(x)$ is integrable over $(0,+\infty)$. We have:
\begin{eqnarray*} \log(\vartheta_{\Delta}(\lambda))&=& -\Delta\int_0^{+\infty}\frac{1-e^{-\lambda x}}x xn_Z(x)dx=
-\Delta\int_0^{+\infty}(\int_0^{\lambda}e^{-sx}ds) xn_Z(x)dx\\ &=& -\Delta \int_0^{\lambda}\left( 
\int_0^{+\infty}e^{-sx}xn_Z(x)dx\right)ds. \end{eqnarray*}
Hence,
$$\psi_{\Delta}(u)=\exp\left(  -\Delta \int_0^{u^2/2}  \left( \int_0^{+\infty}
    e^{-sx}g_Z(x)dx\right) ds\right).$$
Moreover, it is possible to relate the L\'evy density $n_{L}$ of $(L_t)$ with the
L\'evy density $n_{Z}$ of $(Z_t)$ as follows. Consider $f$ a non negative
function on $\mathbb R$, with $f(0)=0$. Given the whole path $(Z_t)$, the jumps $\delta
L_s=W_{Z_s}-W_{Z_{s-}}$   are   centered   Gaussian  with   variance   $\delta
Z_s$. Hence, 
\begin{eqnarray*}
{\mathbb    E}(\sum_{s\le    t}f(\delta L_s))&=   &\sum_{s\le    t}{\mathbb
  E}(\int_{\mathbb R} f(u)
\exp{(-u^{2}/2 \delta Z_s)}\frac{du}{\sqrt{2 \pi \delta Z_s}})\\
&=& t \int_{\mathbb R} f(u)du \left( \int_0^{+\infty}\exp{(-u^{2}/2 x)}\frac{n_{Z}(x)dx}{\sqrt{2 \pi x}}) \right).
\end{eqnarray*}
This gives $n_{L}(u)= \int_0^{+\infty}\exp{(-u^{2}/2 x)}\frac{n_{Z}(x)dx}{\sqrt{2 \pi x}}$. By the same tools, we see that
$$
{\mathbb  E}(\sum_{s\le  t}|\delta  L_s|)= \sqrt{2/\pi}{\mathbb  E}(\sum_{s\le
  t}\sqrt{\delta Z_s})= t\int_0^{+\infty}\sqrt{x}n_Z(x) dx.
$$
Therefore,  if the above  integral is  finite, the  process $(L_t)$  has finite
variation on compact sets and it holds that $\int_{\mathbb R} |u|n_L(u) du <\infty$.

With $(Z_t)$  a L\'evy-Gamma process, $g_Z(x)=\beta e^{-\alpha x}\1_{x>0}$. Then $\int_0^{+\infty} 
e^{-sx}\beta e^{-\alpha x}dx= \beta/(\alpha+s)$, and $$\psi_{\Delta}(u)=\left(\frac \alpha{\alpha 
+\frac{u^2}2}\right)^{\Delta\beta}.$$
This  model is  the Variance  Gamma stochastic  volatility model  described by
Madan and Seneta~(1990). As noted in K\"{u}chler and
Tappe~(2008), the  Variance Gamma distributions are special  cases of bilateral
Gamma distributions. The condition
$\int_0^{+\infty}\sqrt{x}n_Z(x) dx<\infty$ holds. We can compute, for instance
using the
norming constant for an inverse Gaussian density,
$$
n_L(u)=   \int_0^{+\infty}\exp{(-\frac{1}{2}(\frac{u^{2}}{x}    +   2   \alpha
  x)}\frac{\beta x^{-3/2}dx}{\sqrt{2 \pi }}= \beta (2\alpha)^{1/4})|u|^{-1} \exp{(-(2\alpha)^{1/2}|u|)}
$$

\section{Statistical strategy}

\subsection{Notations}
Subsequently we denote by
$u^*$ the Fourier transform of the function $u$ defined as
$u^*(y)=\int e^{iyx}u(x)dx,$ and by $\|u\|$, $<u,v>$, $u*v$ the quantities
$$\|u\|^2=\int |u(x)|^2 dx,$$ 
$$ <u,v>= \int u(x)\overline{v}(x)dx\mbox{ with }
z\overline{z}=|z|^2 \mbox{ and } u\star v(x)=\int u(y)\bar v(x-y)dy.$$
Moreover, we recall that for any integrable and square-integrable
functions $u, u_1, u_2$,
\begin{equation}\label{Fou1}
(u^*)^*(x)=2\pi u(-x) \mbox{ and }  \langle u_1, u_2\rangle = (2\pi)^{-1} \langle u_1^*,
u_2^*\rangle.
\end{equation}


\subsection{The projection spaces}
As we use projection estimators, we describe now the projection spaces.
Let us define $$\varphi(x)=\frac{\sin(\pi x)}{\pi x} \; \mbox{ and } \;
\varphi_{m,j}(x) = \sqrt{m} \varphi(mx-j),$$ where $m$ is an integer, that can be taken equal to
$2^{\ell}$. It is well known (see Meyer (1990), p.22) that
$\{\varphi_{m,j}\}_{j \in \mathbb{Z}}$ is an orthonormal basis of
the space of square integrable functions having Fourier transforms
with compact support included into $[-\pi m, \pi m]$. Indeed an elementary computation yields 
\begin{equation}\label{phistar}
\varphi^*_{m,j}(x)= \frac{e^{ixj/m}}{\sqrt{m}} \1_{[-\pi m, \pi m]}(x).
\end{equation}
 We denote by
$S_m$ such a space:
\begin{eqnarray*}
S_m&=& {\rm Span}\{\varphi_{_{m,j}}, \; j\in \mathbb{Z}\}\;
=\{h\in \mathbb{L}^2(\mathbb{R}),  \mbox{supp}(h^*)
\subset[-m\pi,m\pi ]\}.
\end{eqnarray*}
We denote by $(S_m)_{m\in \mathcal M_n}$ the collection of linear spaces, where
$${\mathcal M}_n=\{1, \dots, m_n\}$$ and  $m_n\leq n$ is the maximal admissible value of $m$, subject to constraints to be precised later.

In practice, we should consider the truncated spaces $S_m^{(n)} =
{\rm Span}\{\varphi_{_{m,j}}, \; \; j\in \mathbb{Z}, \;\; |j|\leq
K_n \},$ where $K_n$ is an integer depending on $n$, and the
associated estimators. Under  assumption (H6), it is possible and does not change the main part of the study (see Comte {\it et al.}~(2006)). For the  sake of simplicity, we consider here sums over ${\mathbb Z}$.

\subsection{Estimation strategy}

We want to estimate $g$ such that \begin{equation}\label{fondam} 
g^*(x)=-i\frac{\psi'_{\Delta}(x)}{\Delta\psi_{\Delta}(x)}=\frac{\theta_{\Delta}(x)}{\Delta\psi_{\Delta}(x)},\end{equation} 
with $$ \psi_{\Delta}(x)={\mathbb E}(e^{ixZ_1^{\Delta}}),\;\; \theta_{\Delta}(x)=-i\psi'_{\Delta}(x)={\mathbb E}(Z_1^{\Delta}e^{ixZ_1^{\Delta}}).$$

The orthogonal projection $g_m$ of $g$ on
$S_m$ is given by \begin{equation}\label{projg} g_m=\sum_{j\in {\mathbb Z}}
a_{m,j}(g) \varphi_{m,j} \mbox{ with } a_{m,j}(g) =
\int_{\mathbb R} \varphi_{m,j}(x)g(x)dx=\langle \varphi_{m,j},g\rangle.
\end{equation}

We have at hand the empirical versions of $\psi_{\Delta}$ and
$\theta_{\Delta}$: $$\hat\psi_{\Delta}(x)=\frac 1n \sum_{k =1}^n e^{ixZ_{k}^{\Delta}}, \;\;
 \hat \theta_{\Delta}(x)=\frac 1n \sum_{k=1}^n Z_k^{\Delta}e^{ixZ_k^{\Delta}}.$$

Following Neumann~(1997) and Neumann and Reiss~(2007), we truncate $1/\hat\psi_{\Delta}$ and set \begin{equation}\label{psitilde} \frac 1{\tilde
\psi_{\Delta}(x)} = \frac 1{\hat\psi_{\Delta}(x)} \1_{|\hat\psi_{\Delta}(x)|>\kappa_\psi
n^{-1/2}}.\end{equation}

Now, for $t$ belonging to a space $S_m$ of the collection $(S_m)_{m\in {\mathcal M}_n}$, let us define
\begin{equation}\label{critere} \gamma_n(t)=\frac 1n\sum_{k=1}^n\left(\|t\|^2- \frac 1{\pi\Delta}
Z_k^{\Delta}\int e^{ixZ_k^{\Delta}}\frac{t^*(-x)}{\tilde \psi_{\Delta}(x)}dx \right),
\end{equation}
Consider $\gamma_n(t)$ as an approximation of the theoretical contrast
$$
\gamma_n^{th}(t)=\frac 1n\sum_{k=1}^n\left(\|t\|^2- \frac
1{\pi\Delta} Z_k^{\Delta} \int
e^{ixZ_k^{\Delta}}\frac{t^*(-x)}{\psi_{\Delta}(x)}dx \right),$$
The following sequence of equalities, relying on (\ref{Fou1}), explains the choice of the contrast:
\begin{eqnarray*} {\mathbb E}(\frac 1{2\pi\Delta} Z_k^{\Delta}\int e^{ixZ_k^{\Delta}}\frac{t^*(-x)}{\psi_{\Delta}(x)}dx) &=&
\frac 1{2\pi\Delta} \int \theta_{\Delta}(x)\frac{t^*(-x)}{\psi_{\Delta}(x)} dx
 = \frac 1{2\pi}\langle t^*, g^*\rangle =
 \langle
t,g \rangle.\end{eqnarray*} Therefore, we find that
$\mathbb{E}(\gamma_n^{th}(t))= \|t\|^2-2\langle g, t\rangle = \|t-g\|^2 -\|g \|^2$
is minimal when $t=g$. Thus, we define the estimator belonging to $S_m$ by
\begin{equation}\label{truncsanssel} \hat g_m = {\rm Argmin}_{t\in S_m} \gamma_n(t)\end{equation}
This estimator can also be written \begin{equation}\label{explicit} \hat g_m
 = \sum_{j\in {\mathbb Z}} \hat a_{m,j} \varphi_{m,j}, \mbox{ with } \hat
a_{m,j}= \frac 1{2\pi n\Delta}\sum_{k=1}^n
Z_{k}^{\Delta}\int e^{ixZ_k^{\Delta}} \frac{\varphi_{m,j}^*(-x)}{\tilde \psi_{\Delta}(x)}dx,\end{equation}
or $$\hat a_{m,j}= \frac 1{2 \pi \Delta}\int \hat\theta_{\Delta}(x)\frac{\varphi_{m,j}^*(-x)}{\tilde \psi_{\Delta}(x)}dx.$$

\subsection{Risk bound of the collection of estimators}
First, we recall a key Lemma, borrowed from Neumann~(1997) (see his Lemma 2.1):
\begin{lem}\label{neum}
 It holds that, for any $p\geq 1$,
$${\mathbb E}\left(\left|\frac 1{\tilde\psi_{\Delta}(x)}-\frac 1{\psi_{\Delta}(x)}\right|^{2p}\right) \leq C \left( \frac 1{|\psi_{\Delta}(x)|^{2p}} \wedge
\frac{n^{-p}}{|\psi_{\Delta}(x)|^{4p}}\right),$$ where $1/\tilde\psi_{\Delta}$ is defined by (\ref{psitilde}).
\end{lem}
\noindent Neumann's result is for $p=1$ but the extension to any $p$ is straighforward. See also Neumann and Reiss~(2007).
This lemma allows to prove the following risk bound.
\begin{prop}\label{first}
Under Assumptions {\rm (H1)-(H2)}$(4)${\rm -(H3)}, then for all $m$:
\begin{equation}\label{prelimi} {\mathbb E}(\|g-\hat g_m\|^2) \leq \|g-g_m\|^2 + K\frac{{\mathbb E}^{1/2}[(Z_1^{\Delta})^4] \int_{-\pi m}^{\pi m} dx /|\psi_{\Delta}(x)|^2}{n\Delta^2}.
\end{equation}
where $K$ is a constant.
\end{prop}
\noindent It is worth stressing that (H4) is not required for the above result. Therefore, it holds even for exponential decay of $\psi_{\Delta}$. \\

\noindent {\bf Proof of Proposition \ref{first}}. First with Pythagoras Theorem, we have
\begin{equation}\label{decen2} \|g-\hat g_m\|^2=\|g-g_m\|^2 + \|\hat g_m-g_m\|^2.\end{equation}
Let $$a_{m,j}(g)= \frac 1{2\pi \Delta}\int \theta_{\Delta}(x)\frac{\varphi_{m,j}^*(-x)}{\psi_{\Delta}(x)}dx.$$ Then, using Parseval's formula and (\ref{phistar}), we obtain
$$\|\hat g_m-g_m\|^2 = \sum_{j\in {\mathbb Z}} |\hat a_{m,j}-a_{m,j}(g)|^2= \frac 1{2\pi\Delta^2}\int_{-\pi m}^{ \pi m} \left|\frac{\hat\theta_{\Delta}(x)}{\tilde \psi_{\Delta}(x)} - \frac{\theta_{\Delta}(x)}{\psi_{\Delta}(x)}\right|^2 dx.$$
It follows that
\begin{eqnarray} \nonumber {\mathbb E}(\|\hat g_m-g_m\|^2)&\leq & \frac c{\Delta^2} \left\{ \int_{-\pi m}^{\pi m} {\mathbb E} \left|\hat \theta_{\Delta}(x)\left(\frac 1{\tilde \psi_{\Delta}(x)} -\frac 1{\psi_{\Delta}(x)}\right)\right|^2 dx \right. \\  \nonumber && \left. \hspace{1cm} + \int_{-\pi m}^{\pi m}
\frac{{\mathbb E}|\hat\theta_{\Delta}(x)-\theta_{\Delta}(x)|^2}{|\psi_{\Delta}(x)|^2} dx\right\}\\  \label{repere}
&\leq &
\frac c{\Delta^2} \left\{ \int_{-\pi m}^{\pi m} \left( {\mathbb E}(|\hat\theta_{\Delta}(x)-\theta_{\Delta}(x)|^2
\left|\frac 1{\tilde\psi_{\Delta}(x)} -\frac 1{\psi_{\Delta}(x)}\right|^2\right) dx \right. \\  \nonumber &&\hspace{-0.5cm}  \left. + \int_{-\pi m}^{\pi m} \left(\Delta^2 |g^*(x)\psi_{\Delta}(x)|^2  {\mathbb E}\left(\left|\frac 1{\tilde\psi_{\Delta}(x)} -\frac 1{\psi_{\Delta}(x)}\right|^2\right) +\frac 1n \frac{{\mathbb E}[(Z_1^{\Delta})^2]}{|\psi_{\Delta}(x)|^2} \right)dx\right\}
\end{eqnarray}

The Schwarz Inequality yields
$$ {\mathbb E}\left(|\hat\theta_{\Delta}(x)-\theta_{\Delta}(x)|^2
\left|\frac 1{\tilde\psi_{\Delta}(x)} -\frac 1{\psi_{\Delta}(x)}\right|^2\right)\leq 
 {\mathbb E}^{1/2}(|\hat\theta_{\Delta}(x)-\theta_{\Delta}(x)|^4){\mathbb E}^{1/2}\left(
\left|\frac 1{\tilde\psi_{\Delta}(x)} -\frac 1{\psi_{\Delta}(x)}\right|^4\right).$$
Then, with the Rosenthal inequality ${\mathbb E}(|\hat\theta_{\Delta}(x)-\theta_{\Delta}(x)|^4)\leq c{\mathbb E}[(Z_1^{\Delta})^4]/n^2$ and by using Lemma \ref{neum},
$${\mathbb E}\left(\left|\frac 1{\tilde\psi_{\Delta}(x)} -\frac 1{\psi_{\Delta}(x)}\right|^4\right)\leq \frac{C}{|\psi_{\Delta}(x)|^4}$$ so that
\begin{eqnarray*} &&\int_{-\pi m}^{\pi m} \left( {\mathbb E}^{1/2}(|\hat\theta_{\Delta}(x)-\theta_{\Delta}(x)|^4) {\mathbb E}^{1/2}
\left|\frac 1{\tilde\psi_{\Delta}(x)} -\frac 1{\psi_{\Delta}(x)}\right|^4\right) dx\leq  \frac{c{\mathbb E}^{1/2}[(Z_1^{\Delta})^4]}{n}\int_{-\pi m}^{\pi m} \frac{dx }{|\psi_{\Delta}(x)|^2}.\end{eqnarray*}
For the second term, we use Lemma \ref{neum}, to get
$${\mathbb E}\left(\left|\frac 1{\tilde\psi_{\Delta}(x)} -\frac 1{\psi_{\Delta}(x)}\right|^2\right)\leq \frac{Cn^{-1}}{|\psi_{\Delta}(x)|^4}.$$
We obtain 
\begin{equation}\label{fin} {\mathbb E}(\|\hat g_m-g_m\|^2)\leq  
 \frac c{n\Delta^2} ({\mathbb E}^{1/2}[(Z_1^{\Delta})^4]+\Delta^2\|g\|_1^2+ {\mathbb E}[(Z_1^{\Delta})^2])\int_{-\pi m}^{\pi m} \frac{dx}{|\psi_{\Delta}(x)|^2} dx,
\end{equation}
where $\|g\|_1=\int |g(x)|dx$. 
Therefore, gathering (\ref{decen2}) and (\ref{fin}) implies the result. $\Box$

\begin{rmk}
In papers concerned with deconvolution in presence of unknown error densities, the error characteristic function is estimated using a preliminary and independent set of data. This solution is possible here: we may split the sample and use the first half to obtain a preliminary and independent estimator of $\psi_{\Delta}$, and then estimate $g$ from the second half. This would simplify the above proof, but not the study of the adaptive case.
\end{rmk}

\subsection{Discussion about the rates}
Let us study some examples and use (\ref{prelimi}) to get a relevant choice of $m$.
We have $\|g-g_m\|^2=\int_{|x|\geq \pi m} |g^*(x)|^2dx$.
Suppose that $g$ belongs to the Sobolev class $${\mathcal S}(a, L)= \{f, \int
|f^*(x)|^2(x^2+1)^{a}dx\leq L\}.$$
Then, the bias term satisfies $$\|g-g_m\|^2=O(m^{-2a}).$$
Under (H4), the bound of the variance term satisfies $$\frac{\int_{-\pi m}^{\pi m} dx/|\psi_{\Delta}(x)|^2}{n\Delta}= O\left(\frac{m^{2\beta\Delta+1}}{n\Delta}\right).$$
The optimal choice for $m$ is $O((n\Delta)^{1/(2\beta\Delta + 2a+1})$ and the resulting rate for the risk is $(n\Delta)^{-2a/(2\beta\Delta + 2a+1})$. It is worth noting that the sampling interval $\Delta$ explicitely appears in the exponent of the rate. Therefore, for positive $\beta$, the rate is worse for large $\Delta$ that for small $\Delta$.

\noindent $\bullet$ Let us  consider the example of the compound process. In this case $\beta=0$, the upper bound of the mean integrated squared error is of order $O((n\Delta)^{-2a/(2a +1)})$, if $g$ belongs to the Sobolev
class ${\mathcal S}(a, L)$.
Note that if $g$ is analytic i.e. belongs to a class
$${\mathcal A}(\gamma, Q)=\{ f, \int (e^{\gamma x}+e^{-\gamma x})^2|f^*(x)|^2 dx\leq Q\},$$ then the
risk is of order $O(\ln(n\Delta)/(n\Delta))$ (choose $m=O(\ln(n\Delta))$).\\


\noindent $\bullet$ For the Levy Gamma process, we have a more precise result since we have $$|\psi_{\Delta}(u)|=\frac{\alpha^{\beta \Delta}}{(\alpha^2+u^2)^{{\beta\Delta}/2}}, \;\; g^*(x)=\frac \beta{\alpha-ix}.$$
Therefore $\int_{|x|\geq \pi m} |g^*(x)|^2dx=O(m^{-1})$ and $\int_{[-\pi m, \pi m]} dx/|\psi_{\Delta}(x)|^2 = O(m^{2\beta\Delta+1})$. The resulting rate is of order $(n\Delta)^{-1/(2\beta\Delta+2)}$ for a choice of $m$ of order $O((n\Delta)^{1/(2\beta\Delta+2)})$.\\

\noindent $\bullet$ For the Bilateral Gamma process with $(\beta,\alpha)=(\beta',\alpha')$, we have $$\psi_{\Delta}(u)=\frac{\alpha^{\beta\Delta}}{(\alpha^2+u^2)^{\beta\Delta}}, \;\; g^*(x)=\frac \beta{\alpha^2+x^2}.$$
Therefore $\int_{|x|\geq \pi m} |g^*(x)|^2dx=O(m^{-3})$ and $\int_{[-\pi m, \pi m]} dx/|\psi_{\Delta}(x)|^2 = O(m^{4\beta\Delta+1})$. The resulting rate is of order $(n\Delta)^{-3/(4\beta\Delta+4)}$ for a choice of $m$ of order $O((n\Delta)^{1/(4\beta\Delta+4)})$.\\

These examples illustrate that the relevant choice of $m$ depends on the unknown function, in particular on its smoothness. The model selection procedure proposes a data driven criterion to select $m$.\\

\noindent $\bullet$ Consider now the process described in Section \ref{other}.
In that case, it follows from (\ref{equiv}) that  $\int_{[-\pi m, \pi m]} dx/|\psi_{\Delta}(x)|^2 =O(m^{\delta+1/2}\exp(\kappa m^{1/2-\delta}))$  and $\int_{|x|\geq \pi m} |g^*(x)|^2dx=O(m^{-2\delta})$. In this case, choosing $\kappa m^{1/2-\delta}=\ln(n\Delta)/2$ gives the rate $[\ln(n\Delta)]^{-2\delta}$ which is thus very slow, but known to be optimal in the usual deconvolution setting (see Fan~(1991)). This case is not considered in the following for the adaptative strategy since it does not satisfy (H4).

\subsection{Study of the adaptive estimator}

We have to select an adequate value of $m$. For this, we start by defining the term
\begin{equation}\label{phipsi}
\Phi_{\psi}(m)=\int_{-\pi m}^{\pi m} \frac{dx}{|\psi_{\Delta}(x)|^2},
\end{equation}
and the following theoretical penalty
\begin{equation}\label{penalite} {\rm pen}(m)= \kappa (1+ {\mathbb E}[(Z_1^{\Delta})^2]/\Delta ) \frac{\Phi_{\psi}(m)}{n\Delta}. \end{equation}
We set
$$\hat m=\arg\min_{m\in {\mathcal M}_n} \left\{ \gamma_n(\hat
g_m) + {\rm pen}(m)\right\},$$ and study first the ``risk" of $\hat g_{\hat m}$.

\noindent Moreover we need the following assumption on the collection of models ${\mathcal M}_n=\{1, \dots, m_n\}$, $m_n\leq n$:
\begin{itemize}
\item[(H7)$\;\;\;\;\;$] $\exists \varepsilon, 0<\varepsilon <1, \;\; m_n^{2\beta \Delta}\leq Cn^{1-\varepsilon}$,
\end{itemize}
where $C$ is a fixed constant and $\beta$ is defined by (H4).


For instance, Assumption (H7) is fulfilled if:
\begin{enumerate} \item pen$(m_n)\leq C$. In such a case, we have 
$m_n\leq C (n\Delta)^{1/(2\beta\Delta+1)}$.
\item $\Delta$ is small enough to ensure $2\beta\Delta<1$. In such a case we can take ${\mathcal M}_n=\{1, \dots, n\}$.
\end{enumerate}

\begin{rmk} Assumption (H7) raises a problem since it depends on the unknown $\beta$ and concrete implementation requires the knowledge of $m_n$. It is worth stressing that the analogous difficulty arises in deconvolution with unknown error density (see Comte and Lacour~(2008)). In the compound Poisson model, $\beta=0$ and nothing is needed. Otherwise one should at least know if $\psi_{\Delta}$ is in a class of polynomial decay. The estimator $\hat\psi_{\Delta}$ may be used to that purpose and to provide an estimator of $\beta$ (see e.g. Diggle et 
Hall~(1993)).  
\end{rmk}

Let us define $$\theta_{\Delta}^{(1)}(x)= {\mathbb E}(Z_1^{\Delta}\1_{|Z_1^{\Delta}|\leq k_n\sqrt{\Delta}} e^{ixZ_1^{\Delta}}),\;\;\;\; \theta_{\Delta}^{(2)}(x)= {\mathbb E}(Z_1^{\Delta}\1_{|Z_1^{\Delta}|> k_n\sqrt{\Delta}} e^{ixZ_1^{\Delta}})$$
so that $\theta_{\Delta}=\theta_{\Delta}^{(1)}+\theta_{\Delta}^{(2)}$ and analogously 
$\hat \theta_{\Delta}=\hat \theta_{\Delta}^{(1)}+\hat\theta_{\Delta}^{(2)}$.
For any two functions $t, s$ in $S_m$, the contrast $\gamma_n$ satisfies:
\begin{eqnarray} \nonumber\gamma_n(t)-\gamma_n(s)&=&\|t-g\|^2- \|s-g\|^2 -2\nu_n^{(1)}(t-s) -2\nu_n^{(2)}(t-s)
\\ \label{dec} &&  - 2\sum_{i=1}^4 R_n^{(i)}(t-s),
\end{eqnarray}  with
\begin{eqnarray*} \nu_n^{(1)}(t)&=& \frac 1{2\pi \Delta} \int t^*(-x)  \frac{\hat \theta_{\Delta}^{(1)}(x)-\theta_{\Delta}^{(1)}(x)}{\psi_{\Delta}(x)} dx,\\
\nu_n^{(2)}(t)& =&\frac 1{2\pi \Delta}\int t^*(-x) \frac{\theta_{\Delta}(x)}{[\psi_{\Delta}(x)]^2} (\psi_{\Delta}(x)-\hat \psi_{\Delta}(x))dx,\\
R_n^{(1)}(t)&= &\frac 1{2\pi \Delta}  \int t^*(-x) (\hat\theta_{\Delta}(x)-\theta_{\Delta}(x)) \left( \frac 1{\tilde \psi_{\Delta}(x)}-\frac 1{\psi_{\Delta}(x)}\right) dx\\
R_n^{(2)}(t) &=& \frac 1{2\pi \Delta}\int t^*(-x) \frac{\theta_{\Delta}(x)}{\psi_{\Delta}(x)} (\psi_{\Delta}(x)-\hat \psi_{\Delta}(x))\left( \frac 1{\tilde \psi_{\Delta}(x)}-\frac 1{\psi_{\Delta}(x)}\right) dx,\\
R_n^{(3)}(t)&=&\frac 1{2\pi \Delta} \int t^*(-x)  \frac{\hat \theta_{\Delta}^{(2)}(x)-\theta_{\Delta}^{(2)}(x)}{\psi_{\Delta}(x)} dx,\\
R_n^{(4)}(t)&=&-\frac 1{2\pi \Delta}\int t^*(-x) \frac{\theta_{\Delta}(x)}{\psi_{\Delta}(x)} \1_{|\hat \psi_{\Delta}(x)|\leq \kappa_{\psi}/\sqrt{n}} dx.
\end{eqnarray*}

Using this decomposition and Talagrand's inequality, we can prove

\begin{thm}\label{thmpenth}
Assume that assumptions {\rm (H1)-(H2)}$(8)${\rm -(H3)-(H7)} hold.
Then
$${\mathbb E}(\|\hat g_{\hat m}-g\|^2) \leq C \inf_{m\in {\mathcal M}_n} \left( \|g-g_m\|^2 +
{\rm pen}(m)\right) + K\frac{\ln^2(n)}{n\Delta},$$
where $K$ is a constant.
\end{thm}

\begin{rmk}
Assumption (H6) is satisfied for the Levy-Gamma process. For the compound Poisson process, it is equivalent to  $\int x^4f^2(x)dx<+\infty$, where $f$ denotes the density of $Y_i$ (see Section \ref{example}).

\end{rmk}

To get an estimator, we replace the theoretical penalty by: $$\widehat{{\rm pen}}(m)=\kappa'
\left(1+ \frac 1{n\Delta^2} \sum_{i=1}^n
(Z_i^{\Delta})^2\right) \frac{\int_{-\pi m}^{\pi m}dx /|\tilde
\psi_{\Delta}(x)|^2 dx }{n}.$$
In that case we can prove:
\begin{thm}\label{thmpenhat}
Assume that assumptions {\rm (H1)-(H2)}$(8)${\rm-(H3)-(H7)} hold and let $\tilde g =\hat g_{\widehat{\widehat m}}$
be the estimator defined with $\widehat{\widehat m}=\arg\min_{m\in {\mathcal M}_n} (\gamma_n(\hat g_m)+\widehat{{\rm pen}}(m))$.
Then
$${\mathbb E}(\|\tilde g-g\|^2) \leq C \inf_{m\in {\mathcal M}_n} \left( \|g-g_m\|^2 +
{\rm pen}(m)\right) + K'_{\Delta} \frac{\ln^2(n)}{n}$$
where $K'_{\Delta}$ is a constant depending on $\Delta$ (and on fixed quantities but not on $n$).
\end{thm}
Theorem \ref{thmpenhat} shows that the adaptive estimator automatically achieves the best rate that can be hoped. If $g$ belongs to the Sobolev ball ${\mathcal S}(a,L)$, and under (H4), the rate is automatically of order $O((n\Delta)^{-2a/(2\beta \Delta+2a+1)})$. See Section 4.5.

\begin{rmk} \begin{enumerate} \item  It is possible to extend our study of the adaptive estimator to the case $\psi_{\Delta}$ having exponential decay. Note that the faster $|\psi_{\Delta}|$ decays, the more difficult it will be to estimate $g$. 
\item Few results on rates of convergence are available in the literature for this problem. The results of Neumann and Reiss~(2007) are difficult to compare with ours since the point of view is different.
\end{enumerate} \end{rmk}

\section{Proofs}
\subsection{Proof of Theorem \protect \ref{thmpenth}. }

Writing that $\gamma_n(\hat g_{\hat m})+ {\rm pen}(\hat m) \leq
\gamma_n(g_m)+ {\rm pen}(m)$ in view of (\ref{dec}) implies that
\begin{eqnarray*}
\|\hat g_{\hat m}-g\|^2 &\leq & \|g_m-g\|^2 + 2\nu_n^{(1)}(\hat g_{\hat m}-g_m)
 +2\nu_n^{(2)}(\hat g_{\hat m}-g_m) + 2\sum_{i=1}^4 R_n^{(i)}(\hat g_{\hat m}-g_m)\\
  && + {\rm pen}(m)- {\rm pen}(\hat m).
\end{eqnarray*}

Let us take expectations of both sides and bound each r.h.s. term.
\begin{eqnarray*}
|{\mathbb E}(\nu_n^{(1)}(\hat g_{\hat m}-g_m))| &\leq & \frac
1{16} {\mathbb E}(\|g_m-\hat g_{\hat m}\|^2) + 16 {\mathbb
E}\left[\sup_{t\in S_m+ S_{\hat m}, \|t\|=1}
|\nu_n^{(1)}(t)|^2\right] \\ &\leq & \frac 18 {\mathbb E}(\|g-\hat
g_{\hat m}\|^2) +\frac 18 \|g- g_m\|^2\\ && + 16 {\mathbb
E}\left(\sup_{t\in S_{m\vee\hat m}, \|t\|=1}
|\nu_n^{(1)}(t)|^2-p_1(m,\hat m)\right)_+ +16 {\mathbb E}(p_1(m,\hat
m)).
\end{eqnarray*}

The same kind of bounds are obtained for $\nu_n^{(2)}$ and the residuals leading to
\begin{eqnarray}\nonumber
\frac 28 {\mathbb E}( \|\hat g_{\hat m}-g\|^2) &\leq & \frac{14}8 \|g-g_m\|^2 +
16 \sum_{m'\in {\mathcal M}_n} {\mathbb E}\left(\sup_{t\in S_{m\vee m'}, \|t\|=1} |\nu_n^{(1)}(t)|^2-p_1(m,m')\right)_+ \\ \nonumber && +
16 {\mathbb E}\left(\sup_{t\in S_{m\vee\hat m}, \|t\|=1} |\nu_n^{(2)}(t)|^2-p_2(m,\hat m)\right)_+
\\  \nonumber   && + 16\sum_{i=1}^2 {\mathbb E}\left(\sup_{t\in S_{m\vee \hat m}, \|t\|=1}|R_n^{(i)}(t)|^2-p_1(m,\hat m)\right)
\\  \nonumber  && + 16\sum_{i=3}^4 {\mathbb E}\left(\sup_{t\in S_{m_n}, \|t\|=1}|R_n^{(i)}(t)|^2\right) \\
&& \label{etape1} +  {\rm pen}(m) + {\mathbb E}(48p_1(m,\hat m)+
16p_2(m,\hat m) - {\rm pen}(\hat m)).
\end{eqnarray}
Next, definition of pen$(.)$ comes from the following constraint:
\begin{equation}\label{penpm} 48 p_1(m,m')+16p_2(m,m')\leq {\rm pen}(m')+ {\rm
pen}(m).\end{equation} This leads to
$${\rm pen}(m) + {\mathbb E}(48p_1(m,\hat m)+ 16p_2(m,\hat m) - {\rm pen}(\hat m)) \leq 2{\rm pen}(m).$$
First, we apply Talagrand's Inequality recalled in Lemma \ref{Concent} to 
prove the following result:
\begin{prop}\label{nun1} Under the assumptions of Theorem \ref{thmpenth}, define $$p_1(m,m')=(4{\mathbb E}[(Z_1^{\Delta})^2] \int_{-\pi (m\vee m')}^{\pi(m\vee m')} |\psi_{\Delta}(x)|^{-2}dx)/(\pi n \Delta^2),$$ then
\begin{equation}\label{nunu} \sum_{m'\in {\mathcal M}_n} {\mathbb E}\left(\sup_{t\in S_{m\vee m'}, \|t\|=1} |\nu_n^{(1)}(t)|^2-p_1(m,m')\right)_+\leq \frac cn.\end{equation}
\end{prop}
\noindent Next we prove:
\begin{prop}\label{nun2} Under the assumptions of Theorem \ref{thmpenth}, define $p_2(m,m')=0$ if $-a+\beta\Delta\leq 0$ and $p_2(m,m')=(\int_{-\pi (m\vee m')}^{\pi(m\vee m')} |\psi_{\Delta}(x)|^{-2}dx)/n$ otherwise. Then
\begin{equation}\label{nunu2} {\mathbb E}\left(\sup_{t\in S_{m\vee \hat m}, \|t\|=1} |\nu_n^{(2)}(t)|^2-p_2(m,\hat m)\right)_+\leq \frac cn.\end{equation}
\end{prop}
For the residual terms, two type of results can be obtained.
\begin{prop}\label{rn12}
Under the assumptions of Theorem \ref{thmpenth}, for $i=1,2$,
$${\mathbb E}\left(\sup_{t\in S_{m\vee \hat m}, \|t\|=1}[R_n^{(i)}(t)]^2-p_1(m,\hat m)\right)\leq \frac{C}{n\Delta}.$$
\end{prop}
and 
\begin{prop}\label{rn34}
Under the assumptions of theorem \ref{thmpenth}, for $i=3,4$
$${\mathbb E}\left(\sup_{t\in S_{m_n}, \|t\|=1}[R_n^{(i)}(t)]^2\right)\leq c\frac{\ln^2(n)}{n\Delta}.$$
\end{prop}

Then the choice ${\rm pen}(m)$ given by (\ref{penalite}) gives,
following (\ref{etape1}) and (\ref{penpm}),
$$
\frac 14 {\mathbb E}( \|\hat g_{\hat m}-g\|^2) \leq  \frac{7}4 \|g-g_m\|^2 + 2{\rm pen}(m)
+ C\frac{\ln^2(n)}{n\Delta},$$
which is the result. $\Box$

\subsection{Proof of Proposition \ref{nun1}.}
Let $$\omega_t(z)=\frac{z\1_{\{|z|\leq k_n\sqrt{\Delta}\}}}{2\pi \Delta} \int e^{izx}\frac{t^*(-x)}{\psi_{\Delta}(x)}dx$$
and notice that $$\nu_n^{(1)}(t)= \frac 1{n} \sum_{k=1}^n \left[\omega_t(Z_k^{\Delta})- {\mathbb E}(\omega_t(Z_k^{\Delta}))\right].$$
To apply Lemma \ref{Concent}, we compute $M_1, H_1$ and $v_1$ defined therein.
First, we have
\begin{eqnarray*}
{\mathbb E}\left(\sup_{t\in S_{m}, \|t\|=1} |\nu_n^{(1)}(t)|^2\right) &\leq &
 {\mathbb E} \left(\sum_{j\in {\mathbb Z}} |\nu_n^{(1)}(\varphi_{m,j})|^2\right)\\ & = &   {\mathbb E} \left(\frac 1{2\pi \Delta^2} \int_{-\pi m}^{\pi m} \left|\frac{\hat\theta_{\Delta}^{(1)}(x)-
\theta_{\Delta}^{(1)}(x)}{\psi_{\Delta}(x)}\right|^2 dx\right)
\\ &\leq &  \frac{{\mathbb E}[(Z_1^{\Delta})^2]}{2\pi n\Delta^2} \Phi_{\psi}(m),
\end{eqnarray*}
where $\Phi_{\psi}(m)$ is defined in (\ref{phipsi}). 
We can take, for $m^\star=m\vee m'$, $$H_1^2=
\frac{{\mathbb E}[(Z_1^{\Delta})^2]}{2\pi n\Delta^2} \Phi_{\psi}(m^\star).$$

Then it is easy to see that if $\|t\|=1$ and $t\in S_{m^{\star}}$, then
$$|\omega_t(z)|\leq \frac{k_n}{2\pi \sqrt{\Delta}}\int \left|\frac{t^*(-x)}{\psi_{\Delta}(x)} \right|dx \leq
\frac{k_n}{2\pi \sqrt{\Delta}}\sqrt{\Phi_{\psi}(m^{\star})}:=M_1.$$
Lastly, for $t\in S_m, \|t\|=1$, $t=\sum_{j\in {\mathbb Z}} t_{m,j}\varphi_{m,j}$
\begin{eqnarray*}
&&{\rm Var}(\omega_t ( Z_1^{\Delta} )) \leq \frac 1{(2\pi)^2\Delta^2}\iint {\mathbb E}\left( e^{i(u-v)Z_1^{\Delta}} (Z_1^{\Delta}\1_{|Z_1^{\Delta}|\leq k_n\sqrt{\Delta}})^2\right) \frac{t^*(-u)t^*(v)}{\psi_{\Delta}(u)\psi_{\Delta}(-v)} dudv\\ &=&
\frac 1{(2\pi \Delta)^2}\sum_{j,k} t_{m,j}t_{m,k}\iint {\mathbb E}\left( e^{i(u-v)Z_1^{\Delta}} (Z_1^{\Delta}\1_{|Z_1^{\Delta}|\leq k_n\sqrt{\Delta}})^2\right) \frac{\varphi_{m,j}^*(-u)\varphi_{m,k}^*(v)}{\psi_{\Delta}(u)\psi_{\Delta}(-v)} dudv.
\end{eqnarray*}
Denoting by
\begin{equation}\label{hdelta} h_{\Delta}^*(u)={\mathbb E}[e^{iuZ_1^{\Delta}} (Z_1^{\Delta}\1_{|Z_1^{\Delta}|\leq k_n\sqrt{\Delta}})^2], \end{equation}
we obtain:
\begin{eqnarray*}
{\rm Var}(\omega_t ( Z_1^{\Delta} )) \leq
&\leq &
\frac 1{(2\pi \Delta)^2}\left(\sum_{j,k} \left|\iint h_{\Delta}^*(u-v) \frac{\varphi_{m,j}^*(-u)\varphi_{m,k}^*(v)}{\psi_{\Delta}(u)\psi_{\Delta}(-v)} dudv\right|^2\right)^{1/2}\\
&=&
\frac 1{2\pi \Delta^2}\left(\iint_{[-\pi m,\pi m]^2}  \left| \frac{h_{\Delta}^*(u-v)}{\psi_{\Delta}(u)\psi_{\Delta}(-v)}\right|^2 dudv\right)^{1/2}
\end{eqnarray*}
where the last equality follows from the Parseval equality.  Next with the Schwarz inequality and the Fubini theorem, we obtain
\begin{eqnarray*}
{\rm Var}(\omega_t ( Z_1^{\Delta} )) &\leq & \frac 1{2\pi\Delta^2}\left(\iint_{[-\pi m,\pi m]^2}   \frac{|h_{\Delta}^*(u-v)|^2}{|\psi_{\Delta}(u)|^4} dudv\right)^{1/2}
\\  &=& \frac 1{2\pi \Delta^2}\left(\int_{[-\pi m,\pi m]}  \left| \frac{du}{\psi_{\Delta}(u)}\right|^4 du\int |h_{\Delta}^*(z)|^2dz \right)^{1/2}
\\ &\leq & \frac{\sqrt{\int_{-\pi m}^{\pi m} dx/|\psi_{\Delta}(x)|^4}}{2\pi \Delta} \frac{\|h_{\Delta}^*\|}{\Delta}.\end{eqnarray*}
Now we use the following Lemma:
\begin{lem}\label{bornehstar} Under the assumptions of Theorem \ref{thmpenth},
$$\|h^*_{\Delta}\|/\Delta \leq 2\sqrt{\pi} \left(\int x^2g^2(x)dx +
{\mathbb E}[(Z_1^{\Delta})^2]\|g\|^2\right)^{1/2}:=\xi.$$
\end{lem}
Thus, under (H5), $\xi$ is finite.
We set $$v_1=  \frac{\xi\sqrt{\int_{-\pi m^\star}^{\pi m^\star} dx/|\psi_{\Delta}(x)|^4}}{2\pi\Delta} .$$
Therefore, setting $\epsilon^2=1/2$, $$p_1(m,m')= 4{\mathbb E}[(Z_1^{\Delta})^2 /\Delta]\frac{\Phi_{\psi}(m^\star)}{2\pi n\Delta}(=2(1+2\epsilon^2)H_1^2).$$
Using (H4) and the fact that ${\mathbb E}[(Z_1^{\Delta})^2/\Delta]$ is bounded, we find
\begin{eqnarray*}{\mathbb E}\left(\sup_{t\in S_{m^{\star}}, \|t\|=1} |\nu_n^{(1)}(t)|^2-p_1(m,m')\right)_+ &\leq &
C\left(\frac{(m^\star)^{2\beta\Delta+1/2}}{n\Delta} e^{-K\sqrt{m^\star}} \right. \\ && \hspace{2cm} \left.
+ \frac{k_n^2 \Phi_{\psi}(m^\star)}{n^2\Delta} e^{-K'\sqrt{n}/k_n}\right).
\end{eqnarray*}
Here $K=K(c_\psi, C_\psi)$. Moreover, we take 
\begin{equation}\label{choixkn} k_n=K'\sqrt{n}/((2\beta\Delta+3)\ln(n))
\end{equation}  and we obtain
$$\sum_{m'\in {\mathcal M}_n}{\mathbb E}\left(\sup_{t\in S_{m^\star}, \|t\|=1} [\nu_n^{(1)}(t)]^2-p_1(m,m')\right)_+ \leq \frac{K"}{n\Delta}.$$

\subsection{Proof of Proposition \ref{nun2}.}
The study of $\nu_n^{(2)}$ is slightly different. 
\begin{eqnarray*}
{\mathbb E}\left(\sup_{t\in S_{m}, \|t\|=1} |\nu_n^{(2)}(t)|^2\right) &\leq &
\frac{1}{2 \pi n\Delta^2}\int_{-\pi m}^{\pi m}\frac{|\theta_{\Delta}(x)|^2}{|\psi_{\Delta}(x)|^4}dx
= \frac{1}{2 \pi n}\int_{-\pi m}^{\pi m}\frac{|g^*(x)|^2}{|\psi_{\Delta}(x)|^2}dx.
\end{eqnarray*}
With assumptions (H4) and (H5), we can see that if $-a+\beta\Delta\leq 0$, then
$$\int_{-\pi m}^{\pi m}\frac{|g^*(x)|^2}{|\psi_{\Delta}(x)|^2}dx\leq \int_{-\pi m}^{\pi m} |g^*(x)|^2(1+x^2)^a
\frac{(1+x^2)^{-a+\beta\Delta}}{c_{\psi}^2}dx\leq \frac 1{c_\psi^2} \int |g^*(x)|^2(1+x^2)^a dx\leq \frac{L}{c_\psi^2}.$$
In that case, we simply take $p_2(m,m')=0$ and write
$${\mathbb E}\left(\sup_{t\in S_{m\vee \hat m}, \|t\|=1} [\nu_n^{(2)}]^2(t)\right) \leq {\mathbb E}\left(\sup_{t\in S_{m_n}, \|t\|=1} [\nu_n^{(2)}]^2(t)\right) \leq \frac{L}{n c_\psi^2}.$$
Now we study the case $-a+\beta\Delta>0$ and find the constants $H=H_2, v=v_2, \epsilon=\epsilon_2$ to apply Lemma  \ref{Concent}. Consider $$\tilde\omega_t(z)=(1/2\pi\Delta)\int e^{izu}t^*(-u)\{\theta_{\Delta}(u)/[\psi_{\Delta}(u)]^2\} du.$$ As
$$\int_{-\pi m}^{\pi m}\frac{|g^*(x)|^2}{|\psi_{\Delta}(x)|^2}dx\leq \frac{L}{c_\psi^2} m^{-2a+2\beta\Delta},$$
we take $$H_2^2= \frac{L}{2\pi c_\psi^2} \frac{(m^{\star})^{-2a+2\beta\Delta}}n.$$
Next, we have $$M_2=\sqrt{n}H_2$$ and we use the rough bound $v_2=nH_2^2$. Moreover, we take $\epsilon^2_2=(-2a+2\beta\Delta+2)\ln(m^{\star})/K_1$.
There exists $m_0$, such that for $m^\star\geq m_0$, $$2(1+2\epsilon_2^2)H_2^2\leq \Phi_{\psi}(m^\star)/n.$$
We set $p_2(m,m')=\Phi_{\psi}(m^\star)/n.$ Introducing
$$W_n(m,m')=\left[\sup_{t\in S_{m\vee m'}, \|t\|=1} |\nu_n^{(2)}|^2(t)-p_2(m,m')\right]_+,$$
we find that
\begin{eqnarray*}
\sum_{m'\in {\mathcal M}_n} {\mathbb E}(W_n(m,m')) &=& \sum_{m'| m^\star\leq m_0} {\mathbb E}(W_n(m,m'))
+\sum_{m'| m^\star> m_0} {\mathbb E}(W_n(m,m'))\\
&\leq & \sum_{m'| m^\star \leq m_0} [{\mathbb E}(\sup_{t\in S_{m^\star}, \|t\|=1} |\nu_n^{(2)}(t)|^2-2(1+2\epsilon_2^2) H_2^2]_+)\\ && + \sum_{m'| m^\star \leq m_0} |p_2(m,m')-2(1+2\epsilon_2^2)H_2^2| \\ & & +
 \sum_{m'| m^\star> m_0} {\mathbb E}([\sup_{t\in S_{m^\star}, \|t\|=1} |\nu_n^{(2)}(t)|^2-2(1+2\epsilon_2^2) H_2^2
]_+).\end{eqnarray*}
Therefore
\begin{eqnarray*}
\sum_{m'\in {\mathcal M}_n} {\mathbb E}(W_n(m,m')) & \leq &
2\sum_{m'\in {\mathcal M}_n}{\mathbb E}([\sup_{t\in S_{m^\star}, \|t\|=1} |\nu_n^{(2)}(t)|^2-2(1+2\epsilon_2^2) H_2^2]_+)
\\&&  + \sum_{m'| m^\star \leq m_0} |p_2(m,m')-2(1+2\epsilon_2^2)H_2^2| \\
&\leq & 2\sum_{m'\in {\mathcal M}_n}{\mathbb E}([\sup_{t\in S_{m^\star}, \|t\|=1} |\nu_n^{(2)}(t)|^2-2(1+2\epsilon_2^2) H_2^2]_+) + \frac{C(m_0)}n.\end{eqnarray*}
Talagrand's Inequality again can be then applied and gives that
$$\sum_{m'\in {\mathcal M}_n}{\mathbb E}([\sup_{t\in S_{m^\star}, \|t\|=1} |\nu_n^{(2)}(t)|^2-2(1+2\epsilon_2^2) H_2^2]_+) \leq \frac Cn.$$
The result for $\nu_n^{(2)}$ in this case follows then by saying as for $\nu_n^{(1)}$ that
$${\mathbb E}\left(W_n(m,\hat m)\right) \leq
\sum_{m'\in {\mathcal M}_n}{\mathbb E}(W_n(m,m')).$$

\subsection{Proof of Proposition \ref{rn12}.}

First define $\Omega(x)=\Omega_1(x)\cap\Omega_2(x)$ with 
\begin{eqnarray*} \Omega_1(x)&=&\left\{|\hat\theta_{\Delta}(x)-\theta_{\Delta}(x)|\leq 8{\mathbb E}^{1/2}[(Z_1^{\Delta})^2] (\log^{1/2}(n) n^{-1/2}\right\}, \\ \Omega_2(x)&=&\left\{ \left|\frac 1{\tilde \psi_{\Delta}(x)}-\frac 1{\psi_{\Delta}(x)}\right|\leq 1/(\log^{1/2}(n)n^{\omega}|\psi_{\Delta}(x)|^2)\right\}.
\end{eqnarray*} 
Then split: $R_n^{(1)}(t)=R_n^{(1,1)}(t)+R_n^{(1,2)}(t)$ where
$$R_n^{(1,1)}(t)=\frac 1{2\pi \Delta} \int t^*(-x)(\hat\theta_{\Delta}-\theta_{\Delta})(x)\left(\frac 1{\tilde \psi_{\Delta}(x)}-\frac 1{\psi_{\Delta}(x)}\right)\1_{\Omega(x)}dx$$ and $R_n^{(1,2)}(t)$ the integral on the complement of $\Omega(x)$.
$$
{\mathbb E}\left(\sup_{t\in S_{m\vee \hat m}, \|t\|=1}|R_n^{(1)}(t)|^2\right) \leq 
2{\mathbb E}\left(\sup_{t\in S_{m\vee \hat m}, \|t\|=1}|R_n^{(1,1)}(t)|^2\right)+
2 {\mathbb E}\left(\sup_{t\in S_{m_n}, \|t\|=1}|R_n^{(1,2)}(t)|^2\right)
$$

\begin{eqnarray*}
&& {\mathbb E}\left(\sup_{t\in S_{m\vee \hat m},
\|t\|=1}|R_n^{(1,1)}(t)|^2\right)\\ & \leq & \frac 1{2\pi \Delta^2}{\mathbb E}\left(\int_{-\pi (m\vee \hat m)}^{\pi (m\vee \hat m)}  |\hat\theta_{\Delta}(x)-\theta_{\Delta}(x)|^2\left|\frac 1{\tilde\psi_{\Delta}(x)}-\frac 1{\psi_{\Delta}(x)}\right|^2\1_{\Omega(x)} dx \right) \\
 &\leq  &
\frac{8({\mathbb E}[(Z_1^{\Delta})^2]/\Delta)}{2\pi n \Delta}{\mathbb E}\left(\int_{-\pi (m\vee \hat
m)}^{\pi (m\vee \hat m)}  n^{-2\omega}
\frac{dx}{|\psi_{\Delta}(x)|^4} \right) \leq \frac{4{\mathbb E}[(Z_1^{\Delta})^2]}{\pi \Delta}{\mathbb E}(\frac{\Phi_{\psi}(m\vee \hat m)}{n\Delta}) \leq {\mathbb E}(p_1(m,\hat m)),
\end{eqnarray*}
under the condition $-2\omega + (1-\varepsilon)\leq 0$. Therefore we choose $\omega=(1-\varepsilon)/2$. Note that if $\beta=0$ the decomposition is useless and the residual is straightforwardly negligible.

On the other hand, Lemma (\ref{neum}) yields:  $${\mathbb E}^{1/4}\left[\left|\frac 1{\tilde\psi_{\Delta}(x)} - \frac 1{\psi_{\Delta}(x)}\right|^8
\right]\leq \frac{C_{\Delta}}{n|\psi_{\Delta}(x)|^4}.$$
Now, we find
\begin{eqnarray*}
&& {\mathbb E}\left(\sup_{t\in S_{m_n}, \|t\|=1}|R_n^{(1,2)}(t)|^2\right) \\ &\leq&
 \frac 1{2\pi \Delta^2}\int_{-\pi m_n}^{\pi m_n}  {\mathbb P}^{1/2}(\Omega(x)^c){\mathbb E}^{1/4}[(\hat\theta_{\Delta}(x)-\theta_{\Delta}(x))^8] {\mathbb E}^{1/4}\left[\left|\frac 1{\tilde\psi_{\Delta}(x)} - \frac 1{\psi_{\Delta}(x)}\right|^8\right] dx\\ &\leq & \frac{C{\mathbb E}^{1/4}[(Z_1^{\Delta})^8]}{2\pi n^2}\int_{-\pi m_n}^{\pi m_n}\frac{{\mathbb P}^{1/2}(\Omega(x)^c)}{|\psi_{\Delta}(x)|^4}dx
 \\ &\leq & \frac{C{\mathbb E}^{1/4}[(Z_1^{\Delta})^8]n^{2(1-\varepsilon)+1-b}}{n^2}\leq \frac{C'_{\Delta}}{n} \mbox{ if } {\mathbb P}(\Omega(x)^c)\leq n^{-2b} \mbox{ and } 2(1-\varepsilon)-b\leq 0.
 \end{eqnarray*}
We take $b=2(1-\varepsilon)$. 
In fact, 
$${\mathbb P}(\Omega(x)^c)\leq {\mathbb P}(\Omega_1(x)^c) + {\mathbb P}(\Omega_2(x)^c).$$  We use
the Markov Inequality to bound ${\mathbb P}(\Omega_2(x)^c)$:
\begin{eqnarray*}
{\mathbb P}(\Omega_2(x)^c) &\leq & \log^p(n) n^{2p\omega} |\psi_{\Delta}(x)|^{4p}{\mathbb E}
\left(\left|\frac 1{\tilde\psi_{\Delta}(x)} -\frac 1{\psi_{\Delta}(x)}\right|^{2p}\right)
\\ &\leq & \log^p(n) n^{2p\omega-p}.
\end{eqnarray*}
The choice of $p$ is thus constrained by $2p\omega-p=-p(1-2\omega)< -4(1-\varepsilon)$ that is $p>4(1-\varepsilon)/\varepsilon$, e.g. $p=5(1-\varepsilon)/\varepsilon$.

We use the decomposition of $\theta_{\Delta}(x)=\theta_{\Delta}^{(1)}(x)+\theta_{\Delta}^{(2)}(x)$ with 
$$k_n\sqrt{\Delta} = \frac{\sqrt{n{\mathbb E}[(Z_1^{\Delta})^2]}}{8\sqrt{\log(n)}}.$$ We use the Bernstein Inequality to bound ${\mathbb P}(\Omega_1(x)^c)$. If $X_1, \dots, X_n$ are i.i.d. variables with variance less than $v^2$ and such that $|X_i|\leq c$, then for $S_n=\sum_{i=1}^n X_i$, we have:
$${\mathbb P}(|S_n-{\mathbb E}(S_n)|\geq n\epsilon)\leq 2 \exp\left(-\frac{n\epsilon^2/2}{v^2+c\epsilon}\right).$$ This yields
\begin{eqnarray*}
{\mathbb P}(\Omega_1(x)^c) &\leq & {\mathbb P}\left(|\hat\theta^{(1)}_{\Delta}(x)-\theta_{\Delta}^{(1)}(x)|\geq 
4\sqrt{{\mathbb E}[(Z_1^{\Delta})^2]\log(n)/n}\right) \\ && + {\mathbb P}\left(|\hat\theta^{(2)}_{\Delta}(x)-\theta_{\Delta}^{(2)}(x)|\geq 4\sqrt{{\mathbb E}[(Z_1^{\Delta})^2]\log(n)/n}\right)\\ &\leq & n^{-16/3} + \frac{n}{16{\mathbb E}[(Z_1^{\Delta})^2]\log(n)} {\mathbb E}(|\hat\theta^{(2)}_{\Delta}(x)-\theta^{(2)}_{\Delta}(x)|^2)
\\ &\leq & n^{-16/3} + \frac{{\mathbb E}[(Z_1^{\Delta})^2\1_{|Z_1^{\Delta}|\geq k_n\sqrt{\Delta}}]}
{16{\mathbb E}[(Z_1^{\Delta})^2]\log(n)}  \\ &\leq & 
n^{-16/3} + \frac{8^4{\mathbb E}[(Z_1^{\Delta})^6]\log^2(n)}{16{\mathbb E}^3[(Z_1^{\Delta})^2]n^2} \\ 
&\leq & n^{-16/3} + \frac{c}{n^2\Delta^2}.
\end{eqnarray*}
This gives the result of Proposition \ref{rn12} for $R_n^{(1)}$. 
The study of $R_n^{(2)}$ follows the same line and is omitted. 

\subsection{Proof of Proposition \ref{rn34}.}
First we study $R_n^{(3)}$. 
\begin{eqnarray*}
{\mathbb E}\left(\sup_{t\in S_{m_n}, \|t\|=1}|R_n^{(3)}(t)|^2\right) & \leq &
\frac 1{4\pi^2\Delta^2} {\mathbb E}\left[\sup_{t\in S_{m_n}, \|t\|=1}\left|\int (\hat\theta_{\Delta}^{(2)}(x)-\theta_{\Delta}^{(2)}(x)) \frac{t^*(-x)}{\psi_{\Delta}(x)}dx\right|^2\right] \\
&\leq &
\frac 1{2\pi\Delta^2}\int_{-\pi m_n}^{\pi m_n}  {\mathbb E}[|\hat\theta_{\Delta}^{(2)}(x)-\theta_{\Delta}^{(2)}(x)|^2]  \frac{dx}{|\psi_{\Delta}(x)|^2} \\ &=&
\frac 1{2\pi\Delta^2}\int_{-\pi m_n}^{\pi m_n}  \frac{{\rm Var}(Z_1^{\Delta}\1_{|Z_1^{\Delta}|\geq k_n\sqrt{\Delta} })}n \frac{dx}{|\psi_{\Delta}(x)|^2} \\  &\leq &
\frac{{\mathbb E}[(Z_1^{\Delta})^8]\Phi_{\psi}(m_n)}{2\pi n k_n^6 \Delta^4}
\\ & \leq &
\frac{K {\mathbb E}[(Z_1^{\Delta})^8]\ln^6(n)}{n^{2+\varepsilon} \Delta^4},
\end{eqnarray*}
using the choice of $k_n$ given by (\ref{choixkn}). \\
Next,
\begin{eqnarray*}
{\mathbb E}\left(\sup_{t\in S_{m_n}, \|t\|=1}|R_n^{(4)}(t)|^2\right)& \leq & \frac 1{2\pi\Delta} \int_{-\pi m_n}^{\pi m_n} |g^*(x)|^2 {\mathbb P}(|\hat\psi_{\Delta}(x)|\leq \kappa_{\psi} / \sqrt{n}) dx
\leq \frac c{n\Delta}.
\end{eqnarray*}
If $|\psi_{\Delta}(x)|\geq 2\kappa_{\psi}/\sqrt{n}$, then
\begin{eqnarray*} {\mathbb P}(|\hat \psi_{\Delta}(u)|\leq \kappa_{\psi}n^{-1/2})&\leq & 
{\mathbb P}(|\hat \psi_{\Delta}(u)-\psi_{\Delta}(u)|\leq |\psi_{\Delta}(u)|-\kappa_{\psi}n^{-1/2})\\
&\leq & {\mathbb P}(|\hat \psi_{\Delta}(u)-\psi_{\Delta}(u)|\geq \frac 12 |\psi_{\Delta}(u)|)\\
&\leq & \exp(-cn|\psi_{\Delta}(u)|^2)
\end{eqnarray*} for some $c>0$, where the last inequality follows from Bernstein's Inequality.

Now, it follows from (H4) that $|\psi_{\Delta}(u)|\geq c_{\psi}(1+u^2)^{-\Delta \beta/2}$. Therefore, for $|u|\leq \pi m_n$ with $m_n^{2\beta\Delta} \leq Cn^{1-\varepsilon}$ by (H7), $$|\psi_{\Delta}(u)|\geq c'm_n^{-\beta\Delta}\geq 2\kappa_{\psi}n^{-1/2}.$$

Moreover, with the previous remarks, $\exp(-cn|\psi_{\Delta}(u)|^2)\leq \exp(-cn^{\varepsilon})$ and thus
$$\int_{-\pi m_n}^{\pi m_n} |g^*(x)|^2 {\mathbb P}(|\hat\psi_{\Delta}(x)|\leq \kappa / \sqrt{n}) dx\leq \|g^*\|^2\exp(-cn^{\varepsilon}).$$
Therefore
$$ {\mathbb E}\left(\sup_{t\in S_{m_n}, \|t\|=1}|R_n^{(4)}(t)|^2\right)\leq \frac c{n\Delta}.$$

\subsection{Proof of Lemma \ref{bornehstar}}
Let us denote by $P_{\Delta}$ the distribution of $Z_1^{\Delta}$ and define
$\mu_{\Delta}(dz)= \Delta^{-1}z P_{\Delta}(dz)$. Let us set
$\mu(dx)=g(x)dx$. Equation (\ref{fondam}) states that
$$
\mu_{\Delta}^{*}= \mu^{*} P_{\Delta}^{*}.
$$
Hence, $\mu_{\Delta}= \mu \star P_{\Delta}$. Therefore,
$\mu_{\Delta}$ has a density given by
$$
 \int g(z-y) P_{\Delta}(dy)= \E g(z- Z_1^{\Delta}).
$$
Moreover, we have, for any compactly supported function  $t$:
\begin{equation*}
\frac{1}{\Delta} \E(Z_1^{\Delta} t(Z_1^{\Delta}))= \int t(z) \E g(z- Z_1^{\Delta}) dz= \int
\E(t(x+Z_1^{\Delta}) g(x) dx.
\end{equation*}
Hence, we apply first Parseval formula:
\begin{eqnarray*} \|h_{\Delta}^*\|^2&=& \int |h_{\Delta}^*(x)|^2dx=2\pi \int
h_{\Delta}^2(x)dx =2\pi \Delta \int z^2\1_{|z|\leq
k_n\sqrt{\Delta}}{\mathbb E}^2(g(z-Z_1^{\Delta})) dz \\ &\leq &
2\pi \Delta {\mathbb E}\left(\int z^2\1_{|z|\leq
k_n\sqrt{\Delta}}g^2(z-Z_1^{\Delta})dz \right)
\\& \leq & 2\pi \Delta {\mathbb E}\left( \int (x+Z_1^{\Delta})^2 g^2(x)dx \right) \leq
4\pi\Delta {\mathbb E}\left(\int (x^2+(Z_1^{\Delta})^2)g^2(z)dz \right) \\ &\leq & 4\pi
\Delta\left(  \int x^2g^2(x)+{\mathbb
E}[(Z_1^{\Delta})^2]\|g\|^2\right).
\end{eqnarray*}
This ends the proof. $\Box$

\subsection{Proof of Theorem \ref{thmpenhat}}

Let us define the sets $$\Omega_1=\left\{ \forall m\in {\mathcal M}_n, \int_{-\pi m}^{\pi m} \left|\frac 1{\tilde\psi_{\Delta}(x)} - \frac 1{\psi_{\Delta}(x)}\right|^2dx \leq k_1 \int_{-\pi m}^{\pi m} \frac {dx}{|\psi_{\Delta}(x)|^2}\right\}$$ and
$$ \Omega_2=\left\{ \left|\frac{\frac 1n \sum_{i=1}^n [Z_i^{\Delta }]^2}{{\mathbb E}[(Z_i^{\Delta})^2]}-1\right|\leq k_2\right\}.$$
Take $0<k_1<1/2$ and $0<k_2<1$. On $\Omega_1$, we have, $\forall m\in {\mathcal M}_n$, $$\int_{-\pi m}^{\pi m} \frac{dx}{|\tilde \psi_{\Delta}(x)|^2} \leq (2k_1+2)\int_{-\pi m}^{\pi m} \frac{dx}{|\psi_{\Delta}(x)|^2} \mbox{ and }
\int_{-\pi m}^{\pi m} \frac{dx}{|\psi_{\Delta}(x)|^2} \leq \frac 2{1-2k_1}\int_{-\pi m}^{\pi m} \frac{dx}{|\tilde \psi_{\Delta}(x)|^2}$$
and on $\Omega_2$, we find $$\frac 1n \sum_{i=1}^n [Z_i^{\Delta }]^2\leq (1+k_2){\mathbb E}[(Z_1^{\Delta})^2] \mbox{ and } {\mathbb E}[(Z_1^{\Delta})^2]\leq \frac 1{1-k_2}\frac 1n \sum_{i=1}^n [Z_i^{\Delta }]^2.$$
Il follows that, on $\Omega_1\cap\Omega_2:=\Omega_{1,2}$, we can choose $\kappa'$ large enough to ensure $$48p_1(m,\widehat{\widehat m})+16p_2(m,\widehat{\widehat m})+\widehat{{\rm pen}}(m)-\widehat{{\rm pen}}(\widehat{\widehat m}) \leq C(a,b){\rm pen}(m).$$
This allows to extend the result of Theorem \ref{thmpenth} as follows: $\forall m\in {\mathcal M}_n$, $${\mathbb E}\left( \|\tilde g-g\|^2\1_{\Omega_{1,2}}\right)\leq C\left( \|g-g_m\|^2 + {\rm pen}(m)\right) + \frac{K\ln^2(n)}{n\Delta}.$$
Next we need to prove that \begin{equation}\label{compl} {\mathbb E}\left( \|\tilde g-g\|^2\1_{\Omega_{1,2}^c}\right)\leq \frac{K'}n.\end{equation}
First, we prove that ${\mathbb P}(\Omega_{1,2}^c)\leq c/n^2$ by proving that ${\mathbb P}(\Omega_1^c)\leq c/n^2$ and
${\mathbb P}(\Omega_2^c)\leq c/n$.

\begin{eqnarray*}
{\mathbb P}((\Omega_1)^c) &\leq & \sum_{m\in {\mathcal M}_n} {\mathbb P}\left(
\int_{-\pi m}^{\pi m} \left|\frac 1{\tilde\psi_{\Delta}(x)} - \frac 1{\psi_{\Delta}(x)}\right|^2dx > k_1 \int_{-\pi m}^{\pi m} \frac {dx}{|\psi_{\Delta}(x)|^2}\right) \\ &\leq & \sum_{m\in {\mathcal M}_n} {\mathbb E}\left[\left(\frac{
\int_{-\pi m}^{\pi m} \left|\frac 1{\tilde\psi_{\Delta}(x)} - \frac 1{\psi_{\Delta}(x)}\right|^2dx}{k_1 \Phi_{\psi}(m)}\right)^p\right]\\ &\leq & \sum_{m\in {\mathcal M}_n} \frac{(2\pi m)^{p-1}}{(k_1\Phi_{\psi}(m))^p} {\mathbb E}\left( \int_{-\pi m}^{\pi m} \left| \frac 1{\tilde\psi_{\Delta}(x)} - \frac 1{\psi_{\Delta}(x)}\right|^{2p}dx\right)\\ &\leq &
\sum_{m\in {\mathcal M}_n} C_pm^{p-1}n^{-p}\frac{\int_{-\pi m}^{\pi m} dx/|\psi_{\Delta}(x)|^{4p}}{(\Phi_{\psi}(m))^p} \\ &\leq & \sum_{m\in {\mathcal M}_n} C'_p n^{-p} m^{(p-1) -p(2\beta\Delta+1) + 4p\beta\Delta+1}
= \sum_{m\in {\mathcal M}_n} C'_p m^{2p\beta\Delta}n^{-p}\\ 
&\leq & C" n^{1-p+p(1-\varepsilon)}\leq C"n^{1-p\varepsilon}.
\end{eqnarray*}
As $m^{2\beta\Delta+1}/(n\Delta)$ is bounded $m^{2p\beta\Delta}n^{-p}=O(n^{2p\beta\Delta/(2\beta\Delta+1) -p})=O(n^{-p/(2\beta\Delta+1)})$. Therefore, choosing $p=3/\varepsilon$ ensures that $n^{1-p\varepsilon)}=n^{-2}$ and ${\mathbb P}(\Omega_1^c)\leq C/n^2$.

On the other hand, $${\mathbb P}[\Omega_2^c]\leq \frac 1{k_2^p{\mathbb E}^p[(Z_1^{\Delta})^2]} {\mathbb E}\left(\left|\frac 1n\sum_{i=1}^n [(Z_i^{\Delta})^2-{\mathbb E}[(Z_1^{\Delta})^2]]\right|^p\right).$$
Here the choice $p=4$ gives ${\mathbb P}[\Omega_2^c]=O(1/n^2)$ with a simple variance inequality, provided that ${\mathbb E}[(Z_1^{\Delta})^8]<+\infty$.

\noindent Next, we write that $$\|g-\tilde g\|^2=\|g-g_{\widehat{\widehat m}}\|^2 + \|g_{\widehat{\widehat m}}-\hat g_{\widehat{\widehat m}}\|^2\leq \|g\|^2 + \sum_{j\in {\mathbb Z}} |\hat a_{\widehat{\widehat m},j}- a_{\widehat{\widehat m},j}(g)|^2$$ and
\begin{eqnarray*} \sum_{j\in {\mathbb Z}} |\hat a_{\widehat{\widehat m},j}- a_{\widehat{\widehat m},j}(g)|^2&=&\sum_{j\in {\mathbb Z}} |\nu_n^{(1)}(\varphi_{\widehat{\widehat m},j}) + \nu_n^{(2)}(\varphi_{\widehat{\widehat m},j}) + \sum_{k=1}^4 R_n^{(k)}(\varphi_{\widehat{\widehat m},j})|^2 \\
&\leq & C\sum_{j\in {\mathbb Z}}\{ |\nu_n^{(1)}(\varphi_{\widehat{\widehat m},j})|^2 + |\nu_n^{(2)}(\varphi_{\widehat{\widehat m},j})|^2 + \sum_{k=1}^4 |R_n^{(k)}(\varphi_{\widehat{\widehat m},j})|^2\}\\
&= & C\left\{ \sup_{t\in S_{\widehat{\widehat m}}, \|t\|=1} |\nu_n^{(1)}(t)|^2
+\sup_{t\in S_{\widehat{\widehat m}}, \|t\|=1} |\nu_n^{(2)}(t)|^2 \right.\\
&& \left. + \sum_{k=1}^4 \sup_{t\in S_{\widehat{\widehat m}}, \|t\|=1}
|R_n^{(k)}(t)|^2 \right\}
\end{eqnarray*}
It follows that, ${\mathbb E}(\|g\|^2\1_{\Omega_{1,2}^c})= \|g\|^2{\mathbb P}(\Omega_{1,2}^c)\leq c/n$, and for $k=3,4$,
$${\mathbb E}\left(\sup_{t\in S_{\widehat{\widehat m}}, \|t\|=1} |R_n^{(k)}(t)|^2\1_{\Omega_{1,2}^c}\right) \leq
{\mathbb E}\left(\sup_{t\in S_{m_n}, \|t\|=1} |R_n^{(k)}(t)|^2\right)\leq C/n$$
as it has been proved previously.
Lastly,
\begin{eqnarray*} {\mathbb E}\left(\sup_{t\in S_{\widehat{\widehat m}}, \|t\|=1} |\nu_n^{(1)}(t)|^2\1_{\Omega_{1,2}^c} \right)
&\leq &{\mathbb E}\left(\sup_{t\in S_{\widehat{\widehat m}}, \|t\|=1} \left\{|\nu_n^{(1)}(t)|^2-{\rm pen}(\widehat{\widehat m})\right\}\right)_+  \\
&& + {\mathbb E}\left({\rm pen}(\widehat{\widehat m}) \1_{\Omega_{1,2}^c} \right)\\ &\leq & c( \frac 1{n\Delta} + n{\mathbb P}(\Omega_{1,2}^c))\leq \frac{c'}n
\end{eqnarray*}
using the proof of Theorem \ref{thmpenth} and the fact that pen$(.)$ is less than $O(n)$. The same line can be followed for the other terms.

\section{Appendix}
\begin{lem}
\label{Concent} Let $Y_1, \dots, Y_n$ be independent random
variables, let $\nu_{n,Y}(f)=(1/n)\sum_{i=1}^n [f(Y_i)-{\mathbb E}(f(Y_i))]$ and
let ${\mathcal F}$ be a countable class of
uniformly bounded measurable functions. Then for $\epsilon^2>0$
\begin{eqnarray*}
\mathbb{E}\Big[\sup_{f\in {\mathcal
F}}|\nu_{n,Y}(f)|^2-2(1+2\epsilon^2)H^2\Big]_+ &\leq &  \frac
4{K_1}\left(\frac vn e^{-K_1\epsilon^2 \frac{nH^2}v} +
\frac{98M^2}{K_1n^2C^2(\epsilon^2)} e^{-\frac{2K_1
C(\epsilon^2)\epsilon}{7\sqrt{2}}\frac{nH}{M}}\right),
\end{eqnarray*}
with $C(\epsilon^2)=\sqrt{1+\epsilon^2}-1$, $K_1=1/6$, and
$$\sup_{f\in {\mathcal F}}\|f\|_{\infty}\leq M, \;\;\;\;
\mathbb{E}\Big[\sup_{f\in {\mathcal F}}|\nu_{n,Y}(f)|\Big]\leq
H, \; \sup_{f\in {\mathcal F}}\frac{1}{n}\sum_{k=1}^n{\rm
Var}(f(Y_k)) \leq v.$$\end{lem} This result follows from the
concentration inequality given in Klein and Rio~(2005) and
arguments in Birg\'e and Massart~(1998) (see the proof of their
Corollary 2 page 354). It can be extended to the case where ${\mathcal F}$ is a unit ball of a linear space.


{\small
\begin{thebibliography}{99}
\bibitem{BS} Barndorff-Nielsen O.E. and Shephard N. (2001). Modelling by L\'evy processes for financial econometrics. In: L\'evy processes. Theory and Applications (Barndorff-Nielsen O.E., Mikosch T., Resnick S.L.), 283-318.

\bibitem{BBM} Barron, A.R., Birg\'e, L. and Massart, P.~(1999).  Risk bounds for
model selection via penalization. {\em Probab. Theory Relat.
Relat. Fields.} {\bf 97} 113-150.

\bibitem{BB} Basawa I.V. and Brockwell P.J. (1982). Nonparametric estimation for nondecreasing L\'evy processes. {\em J.R. Statist. Soc., B}, {\bf 44}, 2, 262-269.

\bibitem{Ber} Bertoin J. (1996). L\'evy processes. {\em Cambridge University Press}.

\bibitem{BM} Birg\'e, L. and Massart, P.~(1998).
Minimum contrast estimators on sieves: Exponential bounds and
rates of convergence. {\em Bernoulli} {\bf 4}, 329-375.

\bibitem{CL} Comte F. and Lacour C. (2008). Deconvolution with estimated error. 
\bibitem{CRT} Comte, F., Rozenholc, Y. and  Taupin, M.-L.~(2006) Penalized contrast estimator for adaptive density deconvolution.  {\em Canad. J. Statist.} {\bf 34 }, 431-452.

\bibitem{CT} Cont R. and Tankov P. (2004). Financial modelling with jump processes. In: Financial Mathematics Series, Chapman \& Hall, CRC, Boca Raton.

\bibitem{CS} Cs\"{o}rg\"{o}, S.~(1981) Limit behaviour of the empirical characteristic function.  {\em Ann. Probab.} {\bf 9}, 130-144.
\bibitem{DH} Diggle, P.J. and  Hall, P.~(1993) A Fourier approach to nonparametric deconvolution of a density estimate.  {\em J. Roy. Statist. Soc. Ser. B} {\bf  55}, 523--531. 
\bibitem{Eb}  Eberlein E. and  Keller U.  (1995). Hyperbolic  distributions in
  finance. {\em Bernoulli} {bf 1}, 3, 281-299.
\bibitem{Fer}  Ferguson ,  T.S. and  Klass, M.J.  (1972). A  representation of
  independent increment processes without Gaussian component.
   {\em Ann. Math. Statist.} {\bf 43}, 1634-1643.
\bibitem{FLH} Figueroa-L\'opez J.E. and Houdr\'e C. (2006). Risk bounds for the nonparametric estimation of L\'evy processes. IMS Lecture Notes-Monograph Series High dimensional probability {bf 51}, 96-116.

\bibitem{KR} Klein, T. and  Rio, E.~(2005). Concentration around the mean for
maxima of empirical processes. {\em Ann. Probab.} {\bf 33}
1060-1077.
\bibitem{KT1}  K\"uchler, U. and Tappe, S. (2008)
Bilateral Gamma distributions and processes in financial mathematics.
{\em Stochastic Processes and Their Applications}, {\bf 118}, 261-283.

\bibitem{LT} Ledoux, M. and Talagrand, M. (1991). {\em Probability in
Banach spaces. Isoperimetry and processes.} Springer-Verlag,
Berlin, Ergebnisse der Mathematik und ihrer Grenzgebiete. {\bf 3}.
Folge, 23.
\bibitem{MS} Madan, D.B. and Seneta, E.~(1990) The Variance Gamma (V.G.) Model for Share Market Returns
{\em The Journal of Business} {\bf 63}, No. 4, 511-524.
\bibitem{MEY} Meyer, Y.~(1990) {\em Ondelettes et opérateurs. } I. Hermann, Paris.

\bibitem{NEU} Neumann, M.~(1997) On the effect of estimating the error density in nonparametric deconvolution.  {\em J. Nonparametr. Statist.} {\bf  7}, 307-330.
\bibitem{NR} Neumann, M. and Reiss, M.~(2007). Nonparametric estimation for L\'evy processes from low-frequency observations. {\it Working paper}, ArXiv:0709.2007v1.
\bibitem{Sa} Sato K.I. (1999). {\em L\'evy processes and infinitely divisible distributions.} Cambridge Studies in Advanced Mathematics, 68. Cambridge University Press, Cambridge. 
\bibitem{WK} Watteel R.N. and Kulperger R.J. (2003). Nonparametric estimation of the canonical measure for infinitely divisible distributions. {\em  Journal of Statistical Computation and Simulation}, {\bf 73}, 7, 525-542.
\end{thebibliography}}
\end{document}